\documentclass[12pt]{article}
\usepackage{amsmath,amsfonts,amssymb,amscd,amsthm,amsbsy,epsf,graphics}
\usepackage{amssymb,amsmath,latexsym, amscd, graphicx}
\usepackage{epstopdf}

\newtheorem{theorem}{Theorem}
\newtheorem{proposition}{Proposition}

\newtheorem{corollary}{Corollary}

\def\e{\epsilon}

\def\P{\mathcal P}

\def\R{\mathbb{R}}

\def\Z{\mathbb{Z}}

\def\qed{{\hspace{2mm}{\small $\diamondsuit$}}}

\begin{document}

\title{A topological view of ordered groups}

\author{Dale Rolfsen}

\maketitle

\begin{abstract} In this expository article we use topological ideas, notably compactness, to establish certain basic properties of orderable groups.  Many of the properties we'll discuss are well-known, but I believe some of the proofs are new.
These will be used, in turn, to prove some orderability results, including the left-orderability of the group of PL homeomorphisms of a surface with boundary, which are fixed on at least one boundary component.
\end{abstract}

\section{Orderable groups}
A group $G$ is {\em left-orderable} if there is a strict total ordering $<$ of its elements which is left-invariant, that is $g<h$ implies $fg<fh$ for all $f,g,h \in G$. 

It is easy to check that, given a left-ordering $<$ of $G$, the {\em positive cone} $P = P_< := \{g \in G| 1 < g\}$ satisfies:

(1) $P\cdot P \subset P$ (that is, $P$ is a sub-semigroup)

(2) For each $g \in G$, exactly one of $g=1, g \in P,$ or $g^{-1}\in P$ holds. ($G$ is partitioned: $G = \{1\} \sqcup P \sqcup P^{-1}$)

Conversely, given a subset $P$ of $G$ satisfying (1) and (2), one can define a left-ordering of $G$ by 
$$g< h \iff g^{-1}h \in P.$$  The correspondence $< \to P_<$ is a bijection between the set of left orderings and the set of subsets of $G$ satisfying (1) and (2).  It is sometimes more convenient to consider the left-ordering to be a subset of $G$, in other words an element of the power set $\P(G)$, rather than a relation on the elements of $G$; we will adopt this viewpoint.

It is easy to see that a left-orderable group is also right-orderable; the criterion $g \prec h \iff gh^{-1} \in P$ defines a right-ordering with the same positive cone.  In fact the literature is about evenly divided between discussing left- and right-ordered groups.  If $G$ has a left-ordering which is also right-invariant, we say it is {\em bi-orderable}.  This is equivalent to the positive cone being normal:

(3) $g^{-1}Pg \subset P$ for all $g \in G$.

Useful reference books on orderable groups are \cite{MR77},  \cite{Gla99} and \cite{KM96}.  The article
\cite{Conrad59} is also highly recommended.

\subsection{Algebraic properties of orderable groups}

Knowing that a group is orderable tells us that it has certain special algebraic properties.

$\bullet$ Left-ordered groups $G$ are torsion-free.  

For if $1<g$, then $g<g^2$, $g^2<g^3$ and by transitivity $1<g^n$ for all $n$.  Similarly, if $g<1$ no positive power of $g$ can equal the identity.

$\bullet$ Suppose $\varphi : G \to H$ is a surjective homomorphism with kernel $K$.  If $K$ and $H$ are left-orderable, then so is $G$.

In fact, one can take a positive cone for $G$ the union of the positive cone of $K$ and the preimage under $\varphi$ of the positive cone of $H$.  This does not hold for biorderable groups unless there is a biordering of $K$ invariant under conjugation by elements of $G$.

$\bullet$ Left-orderable groups $G$ satisfy the zero-divisor conjecture, that is, the group ring
$\Z G$ has no zero divisors.

The proof is not difficult, but we omit it here.  It is unknown whether the integral group ring of an arbitrary torsion-free group can have zero divisors.

$\bullet$  If $G$ is left-orderable and $H$ is any group, and $\Z G$ and $\Z H$ are isomorphic as rings, then  $G$ and $H$ are isomorphic as groups.  This is proved in \cite{LR}.

$\bullet$ Bi-ordered groups do not have generalized torsion: if $g$ is not the identity, then any product of conjugates of $g$ cannot be the identity.

This is because if $g>1$ such a product must also be positive, and if $g<1$ the product will be less than the identity too. 

$\bullet$  Bi-ordered groups have unique roots: $g^n = h^n, \; n>0 \implies g=h$.

To see this, one easily checks that in a biordered group inequalities multiply: $g < h$ and $g'<h'$ imply $gg'<hh'$ (this doesn't necessarily hold in a left-ordered group).  So if $g<h$ we conclude $g^2<h^2$, $g^3<h^3$, etc.  The powers can never be equal.

$\bullet$  In a bi-ordered group, if $g^n$ commutes with $h$ for some $n > 0$, then $g$ commutes with $h$.

For if $g$ and $h$ do not commute, say $g < h^{-1}gh$.  Multiply this inequality by itself repeatedly to conclude $g^n<h^{-1}g^nh$.

\subsection{Examples}  Many groups of interest to topologists are orderable.

$\bullet$ $\Z^n$ is bi-orderable, as an additive group.  For example, use the lexicographic ordering.  There are uncountably many possible orderings of $\Z^n$ for $n \ge 2$.  For $\Z^2$, one may take all integral lattice points to one side of a line through the origin with irrational slope as an example of a positive cone.

$\bullet$ Free groups are bi-orderable.  More generally, Vinogradov \cite{Vin} proved the free product of biorderable groups is biorderable.

$\bullet$ Braid groups are left-orderable (Dehornoy \cite{Deh}) but not bi-orderable for more than two strands.

That wonderful and surprising result is what first got me interested in orderable groups. 

$\bullet$ Pure braid groups are bi-orderable. \cite{RZ98} \cite{KR} 

$\bullet$ Fundamental groups of surfaces are bi-orderable, except the Klein bottle group
$\langle x, y| x^{-1}yx = y^{-1} \rangle$
which is only left-orderable, and the projective plane's group which is not even left-orderable, as it is a torsion group.  

This is proved in \cite{RW}.  The Klein bottle group cannot be biordered.  If it were, the defining relation would imply that $y$ is positive if and only if $y^{-1}$ is positive, a contradiction.  However, it is left-orderable, because if one maps it onto $\Z$ by killing the (normal) infinite cyclic subgroup 
$\langle y \rangle$, we have left-orderable kernel and image.

$\bullet$ All classical knot groups are left-orderable and some (but not all) are biorderable.

This is a consequence of a more general result about ordering 3-manifold groups, which we will discuss in Section \ref{3D}.  See also \cite{PR} and \cite{CR}.

$\bullet$ The group $Homeo(I, \partial I)$ of homeomorphisms $h$ of the unit interval $I = [0, 1]$, such that $h(0) = 0$ and $h(1) = 1$, is left-orderable.  Here the group operation is composition.

To see this, choose a well-ordering $r_1 \prec r_2 \prec \cdots$ of the rational numbers in the interval $(0, 1)$.  For two functions
$g, h \in Homeo(I, \partial I)$ declare $g < h$ if and only if $g(r_k) < h(r_k)$ in the usual ordering of $I$, where $r_k$ is the first rational (in the well-ordering) at which the values of $g$ and $h$ differ.  

A similar argument shows that the group of orientation-preserving homeomorphisms of the reals (or the rationals) is left-orderable, that is

$\bullet Homeo^+(\R)$ is left-orderable.

The group $SL(2,\R)$ acts on the circle (for example by fractional linear transformations of 
$\R \cup \{ \infty \}$), and in fact has the homotopy type of $S^1$.  Its universal cover 
$\widetilde{SL}(2, \R)$ is a group which acts on the real numbers by order-preserving homeomorphisms -- it is one of the eight 3-manifold geometries of Thurston \cite{Thu}.  Therefore it may be considered a subgroup of $Homeo^+(\R)$ and we conclude

$\bullet \widetilde{SL}(2, \R)$ is left-orderable.

\section{Topology on the power set}

For any set $X$, one may consider the collection of all its subsets -- that is its power set -- often denoted $\P(X)$ or $2^X$.  This latter notation indicates that the power set may be identified with the set of all functions $X \to \{0, 1\}$ (using von Neumann's  definition $2 :=  \{0, 1\}$), via the characteristic function $\chi_A : X \to  \{0, 1\}$ associated to a subset $A \subset X$ defined by 
$$\chi_A(x) =
  \begin{cases}
     1 \text{ if } x \in A, \\
    0 \text{ if } x \notin A.
   \end{cases}
$$

The set $2^X$ is a special case of a product space: one gives $\{0, 1\}$ the discrete topology, and $2^X$ is considered the product of copies of $\{0, 1\}$ indexed by 
the set $X$.  The product topology is the the smallest topology on the set $2^X$ such that for each $x \in X$ the sets $ \{f \in 2^X : f(x) = 0\}$ and $ \{f \in 2^X : f(x) = 1\}$ are open.  In other notation, the subsets of  $\mathcal{P}(X)$ of the form
$$U_x = \{A \subset X | x \in A\} \quad \text{ and } \quad U_x^c = \{A \subset X | x \notin A\}$$
are open in the ``Tychonoff'' topology on the power set.  Note that the sets $U_x$ and $U_x^c$ are also closed, as they are each other's complement.   A basis for the topology can be gotten by taking finite intersections of various $U_x$ and $U_x^c$.
A famous theorem of Tychonoff asserts that an arbitrary product of compact spaces is again compact.  Since the space $\{ 0, 1 \}$ is compact, we conclude:

\begin{proposition}
The power set $\mathcal{P}(X)$ of any set $X$, with the Tychonoff topology, is compact.  It is also totally disconnected.
\end{proposition}

We recall that a space is said to be {\em totally disconnected} if for each pair of points, there exist disjoint open neighbourhoods of the two points whose union is the whole space.  If $Y$ and $Z$ are distinct points of $P(X)$ (that is, subsets of $X$), take $x\in X$ to be some element of one of those subsets, but not the other; then the sets $U_x$ and $U_x^c$ form such neighbourhoods.
 
If $X$ is finite, then so is $\P(X)$ and the Tychonoff topology is just the discrete topology.  If $X$ is countably infinite, then  $\P(X)$ is homeomorphic to the Cantor space obtained by deleting middle thirds successively of the interval $[0, 1]$.   In particular,
the Tychonoff topology on $\mathcal{P}(X)$ is metrizable when $X$ is countable.

{\bf Example:}
Let $G$ be a group and define $\mathcal{S}(G)$ to be the collection of all sub-semigroups of $G$.  That is, $\mathcal{S}(G) = \{ S \subset G | g, h \in S \implies gh \in S\}$.  Note that 
$\mathcal{S}(G) \subset \mathcal{P}(G)$.  We will argue that 
$\mathcal{S}(G)$ is in fact a {\em closed} subset $\mathcal{P}(G)$.  Consider the complement $\mathcal{P}(G) \setminus \mathcal{S}(G)$.  A subset $Y$ of $X$ belongs
to $\mathcal{P}(G) \setminus \mathcal{S}(G)$ if and only if there exist $g, h \in Y$ with 
$gh \notin Y$.  Therefore
$$\mathcal{P}(G) \setminus \mathcal{S}(G) = 
\bigcup_{g, h \in G}\{U_g \cap U_h \cap U_{gh}^c\}.$$
Each term in the brackets is an open set, by definition, and therefore so is the intersection of the three, and $\mathcal{P}(G) \setminus \mathcal{S}(G)$ is a union of open sets.  It follows that  $\mathcal{S}(G)$ is closed.

\section{The spaces of orderings}

We define the space of left-orderings, $LO(G)$ of a group $G$ to be the collection of 
all subsets $P \subset G$ satisfying (1) and (2) above.  We have just shown that (1) is a closed condition, and a similar argument shows the same for (2).  This proves the following.

\begin{proposition}
$LO(G)$ is a closed subset of $\P(G)$, and is therefore a compact and totally disconnected space (with the subspace topology).
\end{proposition}

This space was introduced in the literature by Adam Sikora \cite{Sik} and has been used to prove some fundamental properties of left-orderable groups: \cite{Mor}, \cite{Lin}.

Although we are considering the topology on left-orderings to be the Tychonoff topology inherited from 
$\P(G)$, there is a natural way to view it in terms of inequalities.
Suppose $<$ is a left-invariant ordering of the group $G$, and suppose we specify a finite number of inequalities $g_1 < h_1,  \dots, g_n < h_n$ which hold.  Then the set of all left-orderings, in which all these inequalities are still true, forms an open neighbourhood of 
$<$ in $LO(G)$.  The set of all such neighbourhoods is a basis for the topology of $LO(G)$.

Similarly, we can define the set $O(G)$ of bi-invariant orderings on the group $G$ to be the collection of subsets $P \subset G$ satisfying (1), (2) and (3) above.  The reader can easily check the following.

\begin{proposition}
$O(G)$ is a closed subset of $LO(G)$, so it is also a compact totally disconnected space. 
\end{proposition}

Of course, for a given $G$ the spaces $LO(G)$ or $O(G)$ may well be empty.

%
%

\section{Testing for orderability}

Suppose we wish to determine if a given group $G$ is left-orderable.  Let's assume for the moment that $G$ is finitely generated, with generators $g_1, \dots, g_n$.  The {\em length} of a group element (relative to the choice of generators) is the smallest integer $k$ such that there is an expression of $g$ in terms of the generators
$$g = g_{i_1}^{\e_1}\cdots g_{i_k}^{\e_k}$$ 
where $\e_i = \pm 1$.  Let $B_k(G)$ denote the set of all elements of $G$ of length at most $k$.  This is a finite set, which includes the identity (length zero) and also is invariant under taking inverses.
It can be regarded as the $k$-ball of the Cayley graph of $G$, relative to the given generators.

Now let us define a subset $Q$ of $B_k(G)$ to be a {\em preorder} of $B_k(G)$ if 

$(1')$ $(Q\cdot Q) \cap B_k(G) \subset Q $ and

$(2')$ $B_k(G) = \{1\} \sqcup Q \sqcup Q^{-1}.$ 

To check whether, for fixed $k$, there exists a preorder of $B_k(G)$ is a finite task.  If one can decide the word problem algorithmically for $G$ (with given generators), then there is an algorithm to decide whether a preorder exists.  Notice that if $P$ is a positive cone 
of a left-ordering of $G$, then $P \cap B_k(G)$ is a preorder of $B_k(G)$, so we conclude the following.

\begin{proposition}
Suppose $G$ is finitely generated by $g_1, \dots, g_n$.   If $G$ is left-orderable, then for every positive integer $k$, $B_k(G)$ admits a preorder.
\end{proposition}

Perhaps surprisingly, there is a converse.

\begin{theorem}\label{partition}
Suppose $G$ is generated by $g_1, \dots, g_n$ and that for all $k \ge 1$, there is a preorder of 
$B_k(G)$.  Then $G$ is left-orderable.
\end{theorem}

We will prove this using compactness of $\mathcal{P}(G)$.  
Consider the set 
$$\mathcal{P}_k = \{ R \subset G |  R \cap B_k(G) \text{ is a preorder of } B_k(G)\}.$$  
One argues as usual that $\mathcal{P}_k$ is a {\em closed} subset of 
$\mathcal{P}(G)$, and by hypothesis $\mathcal{P}_k$ is nonempty.  Note also that 
a preorder of $B_{k+1}(G)$ intersected with $B_k(G)$ becomes a preorder of $B_k(G)$.
That is, we have $\mathcal{P}_{k+1} \subset \mathcal{P}_k$.  Thus the $\mathcal{P}_k$
form a nested descending sequence of nonempty compact subsets of $\mathcal{P}(X)$.
We conclude that 
$$\bigcap_{k=1}^\infty \mathcal{P}_k \ne \emptyset.$$

Also observe that if 
$g, h$ belong to $B_k(G)$ then $gh$ is in $B_{2k}(G)$.  So if 
$P \in \cap_{k=1}^\infty \mathcal{P}_k$ then $P$ is a sub-semigroup.  Similarly $P$ satisfies the partition condition (2) and we conclude that 
$$LO(G) = \bigcap_{k=1}^\infty \mathcal{P}_k \ne \emptyset,$$
completing the proof. \qed

This means that if a finitely-generated group is not left-orderable, then the algorithm described will discover that fact in finite time (although one does not know when!)  Moreover, one can design the algorithm to supply a proof of non-left-orderability if it finds a $B_k(G)$ having no preorder.  On the other hand, if the group under scrutiny {\em is} left-orderable, the algorithm will never end.  An example of such an algorithm, due to Nathan Dunfield, is described in \cite{CD03} and is available from Dunfield's website.  It was used, for example, to find a proof  that the fundamental group of the Weeks manifold -- the smallest volume closed hyperbolic 3-manifold -- is not left-orderable.

The assumption of being finitely-generated is not really essential.

\begin{theorem}\label{fgLO}
A group is left-orderable if and only each of its finitely-generated subgroups is left-orderable.
\end{theorem}

The proof will use the following version of compactness.  A collection of sets is said to have the {\em finite intersection property} if every finite subcollection of the sets has a nonempty intersection.
A space is compact if and only if every collection of closed subsets with the finite intersection property has a nonempty total intersection.

To prove Theorem \ref{fgLO}, consider any finite subset $F$ of the given group $G$ and let
$\langle F \rangle$ denote the subgroup of $G$ generated by $F$.  Define

$$\mathcal{Q}(F) := \{ Q \subset G | Q \cap \langle F \rangle \text{ is a positive cone for }\langle F \rangle \}$$

For each finite $F \subset G$, $\mathcal{Q}(F)$ is a closed subset of $\mathcal{P}(G)$.  The family of all $\mathcal{Q}(F)$, for finite $F \subset G$, is a collection of closed sets which has the finite intersection property, because 
$$ \mathcal{Q}(F_1 \cup F_2 \cup \cdots \cup F_n) \subset \mathcal{Q}(F_1) \cap \mathcal{Q}(F_2) \cap \cdots \cap \mathcal{Q}(F_n).$$

By compactness, $\bigcap_{ \rm{ finite }F \subset G} \mathcal{Q}(F) \ne \emptyset$.

One can easily verify that any element of $\bigcap_{ \rm{ finite }F \subset G} \mathcal{Q}(F)$
is a left-ordering of $G$, completing the proof.  In fact 
$$\bigcap_{ \rm{ finite }F \subset G} \mathcal{Q}(F) = LO(G).$$

\begin{corollary}
An abelian group $G$ is bi-orderable if and only if it is torsion-free.
\end{corollary}

\proof  We need only show that torsion-free abelian groups are left-orderable (which in this case is equivalent to bi-orderable).  But any finitely generated subgroup is isomorphic to 
$\Z^n$ for some $n$, which we have already seen to be bi-orderable.  The result follows from Theorem \ref{fgLO}. \qed

\section{Characterization of left-orderable groups}  

Following \cite{Conrad59}, we have a number of characterizations of left-orderability of a group $G$.  If $X \subset G$, we let $S(X)$ denote the semigroup generated by $X$, that is all elements of $G$ expressible as (nonempty) products of elements of $X$ (no inverses allowed).

\begin{theorem}\label{loequiv}
A group $G$ can be left-ordered if and only if for every finite subset  $ \{x_1, \dots , x_n \}$ of $G$ which does not contain the identity, there exist $\e_i = \pm 1$ such that 
$1 \not\in S(x_1^{\e_1}, \dots , x_n^{\e_n})$.
\end{theorem}

One direction is clear, for if $<$ is a left-ordering of $G$, just choose $\e_i$ so that $x_i^{\e_1}$ is greater than the identity.  For the converse,
by Theorem \ref{fgLO} we may assume that $G$ is finitely generated, and by Theorem \ref{partition} we need only show that each $k$-ball $B_k(G)$, with respect to a fixed finite generating set, has a preorder.  Now consider $ \{x_1, \dots , x_n \}$ to be the entire set $B_k(G) \setminus \{1\}$, and choose $\e_i = \pm 1$ such that 
$1 \not\in S(x_1^{\e_1}, \dots , x_n^{\e_n})$.  

We can easily check that the set $Q := B_k(G) \cap S(x_1^{\e_1}, \dots , x_n^{\e_n})$ is a preorder of $B_k(G)$.  Note that each $x_i$ in the list has its inverse $x_j$ also appearing, and necessarily $\e_i$ and $\e_j$ are opposite in sign, for otherwise 1 would be in the semigroup containing them.  This completes the proof of Theorem \ref{loequiv}.

Another characterization of left-orderability is due to Burns and Hale \cite{BH72}.

\begin{theorem}[Burns-Hale]\label{burnshale}
A group $G$ is left-orderable if and only if for every finitely-generated subgroup 
$H \ne \{ 1 \}$ of $G$, there exists a left-orderable group $L$ and a nontrivial homomorphism $H \to L$.
\end{theorem}

{\bf Proof:}  The forward direction is obvious; just take $L = H$ and use the identity homomorphism.  To prove the other direction, assume the subgroup condition.  According to Theorem \ref{loequiv}, the result will follow if one can show:

Claim:  For every finite subset  $ \{x_1, \dots , x_n \}$ of $G \setminus \{1\}$, there exist $\e_i = \pm 1$ such that  $1 \not\in S(x_1^{\e_1}, \dots , x_n^{\e_n})$.

We will establish this claim by induction on $n$.  It is true for $n = 1$, for $S(x_1)$
cannot contain the identity unles $x_1$ has finite order, which is impossible since the cyclic subgroup 
$\langle x_1\rangle$ must map nontrivially to a left-orderable (hence torsion-free) group.

Next assume the claim true for all finite subsets of $G \setminus \{1\}$ having fewer than $n$ elements, and consider $ \{x_1, \dots , x_n \} \subset G \setminus \{1\}$.  By hypothesis, there is a nontrivial homomorphism 
$$h : \langle x_1, \dots , x_n   \rangle \rightarrow L$$ 
where $(L, \prec)$ is a left-ordered group.  Not all the $x_i$ are in the kernel; we may assume they are numbered so that 
$$h(x_i) 
  \begin{cases}
   \ne 1 \text{ if } i= 1, \dots, r, \\
   = 1  \text{ if } r < i \le n.
   \end{cases}
$$   
Now choose $\e_1, \dots, \e_r$ so that $ 1 \prec h(x_i^{\e_i})$ in $L$ for $i= 1, \dots, r$.
For $i > r$, the induction hypothesis allows us to choose $\e_i = \pm 1$ so that 
$1 \not\in S(x_{r+1}^{\e_{r+1}}, \dots, x_n^{\e_n})$.  We now check that $1 \not\in S(x_1^{\e_1}, \dots , x_n^{\e_n})$ by contradiction.  Suppose that $1$ {\em is} a product of some of the $x_i^{\e_i}$.  If all the $i$ are greater than $r$, this is impossible, as $1 \not\in S(x_{r+1}^{\e_{r+1}}, \dots, x_n^{\e_n})$.  On the other hand if some $i$ is less than or equal to $r$, we see that $h$ must send the product to an element strictly greater than the identity in $L$, again a contradiction.  \qed

A group is said to be {\em indicable} if it has the group of integers $\Z$ as a quotient, and {\it locally indicable} if each of its nontrivial finitely-generated subgroups is indicable.  This notion was introduced by Higman \cite{Higman40} to study zero divisors and units in group rings.

\begin{corollary}
Locally indicable groups are left-orderable.
\end{corollary}

We mention here, without proof, that biorderable groups are locally indicable.  So we have the implications: biorderable $\implies$ locally indicable $\implies$ left-orderable.  Neither of these implications can be reversed.   The braid groups can be used to provide examples.   The 3-strand braid  group $B_3$ is locally indicable but not biorderable, and the 5-strand braid group $B_5$ is left-orderable, but not locally indicable.    In fact the commutator subgroup $[B_5,B_5]$ is finitely generated and perfect \cite{GL}, meaning it equals its own commutator subgroup.  Thus there cannot be a nontrivial homomorphism from $[B_5,B_5]$  to $\Z$, or to any other abelian group.  (See \cite{DDRW} for a more detailed discussion of this.)  Another example of a left-orderable group which is not locally indicable is discussed at the end of Section \ref{3D}.

\begin{corollary}
Suppose $G$ is a group which has a (finite or infinite) family of normal subgroups 
$\{G_\alpha\}$ such that  $\cap_\alpha G_\alpha = \{ 1 \}$.  
If all the factor groups $G/G_\alpha$ are 
left-orderable, then $G$ is left-orderable.
\end{corollary}

This corollary follows, for if $H$ is a finitely generated subgroup of $G$, one can choose 
$\alpha$ for which $H \setminus G_\alpha$ is nonempty.  Then the composition of homomorphisms $H \hookrightarrow G \to G/G_\alpha$ is a nontrivial homomorphism of $H$ to a 
left-orderable group. \qed

\section{Characterization of biorderable groups}

Recall that $P$ is the positive cone of a biordered group $(G,<)$ if and only if it satisfies conditions (1), (2) and (3) cited earlier.  That is, it is a sub-semigroup with the partition property and also normal.
The proof of Theorem \ref{fgLO} adapts easily to a proof of the following.

\begin{theorem}\label{fgBiO}
A group is biorderable if and only each of its finitely-generated subgroups is biorderable.
\end{theorem}

If $G$ is a finitely generated biorderable group, we may consider, as before, the set $B_k(G)$ of all elements of length at most $k$, with respect to some fixed set of generators.  We will define a {\em pre-biorder} of $B_k(G)$ to be a subset $Q$ of $B_k(G)$ satisfying the conditions for a preorder

$(1')$ $(Q\cdot Q) \cap B_k(G) \subset Q $ and

$(2')$ $B_k(G) = \{1\} \sqcup Q \sqcup Q^{-1}.$ 

plus the condition

$(3')$ If $g \in B_k(G)$ then $g^{-1}Qg \cap B_k(G) \subset Q$,

in other words, if one conjugates an element of $Q$ by an element of $B_k(G)$, and the result is still in 
$B_k(G)$, then it must be in $Q$.  Again, checking these conditions, for a fixed $B_k(G)$ is a finite task.
Note that in $(3')$, closure under the other conjugation $gQg^{-1}$ follows, because $B_k(G)$ is closed under taking inverses. 

The following two theorems can be proved in a similar way to their counterparts in the previous section.  We leave the details to the reader. 

\begin{theorem}
A finitely-generated group $G$ is bi-orderable if and only if for every positive integer $k$, the $k$-ball $B_k(G)$ relative to a fixed set of generators admits a pre-biorder.
\end{theorem}

\begin{theorem}
A group $G$ is bi-orderable if and only if for every finite subset 
$\{x_1, \dots, x_n\} \subset G \setminus \{1\}$, there exist $\e_i = \pm 1$
such that $1 \notin S$, where $S$ is the sub-semigroup  generated by the $x_i^{\e_i}$ and their conjugates $x_j^{-1}x_i^{\e_i}x_j$.
\end{theorem}

This is similar to, but sharper than, a characterization due to Fuchs \cite{Fuc}, in which $S$ is replaced by the semigroup generated by the $x_i^{\e_i}$ and their conjugates by {\em all} elements of $G$.   

Note that there is no direct counterpart to the Burns-Hale theorem for biorderable groups.  If there were, then locally indicable groups should be biorderable, which as mentioned above is not always the case.  The inductive step of the proof does not really carry over to the biorderable case, because of all the conjugates which must be considered.

\section{Some applications}

\subsection{3-dimensional manifolds}\label{3D}

Since this is being presented in the memoir of a knot theory conference, it is appropriate to mention the following application of the Burns-Hale theorem, although it already appears in \cite{BRW} and is based on ideas in \cite{HS}.  We outline the proof for the reader's convenience, and since it is a nice application of the Burns-Hale theorem.

\begin{theorem}\label{LO3mfld}
Suppose $M$ is an orientable irreducible 3-manifold.  Then $\pi_1(M)$ is left-orderable if and only if there is a nontrivial homomorphism $h: \pi_1(M) \to L$, where $L$ is a left-orderable group.
\end{theorem}

{\bf Proof:}  The forward direction is obvious.  For the other direction, we will apply the Burns-Hale theorem. If $H$ is a nontrivial finitely-generated subgroup of $\pi_1(M)$, we need to find a nontrivial homomorphism from $H$ to a left-orderable group.
 
Case 1: $H$ has finite index.  This is easy; consider the restriction of $h$ to $H$, which maps $H$ nontrivially to $L$.

Case 2: $H$ has infinite index.   Then there is a covering $p:\tilde{M} \to M$
with $p_*\pi_1(\tilde{M}) = H$.   $\tilde{M}$ is noncompact, but its fundamental group is finitely-generated so, by a theorem of P. Scott \cite{Sco}, there is a compact 3-dimensional submanifold $C \subset \tilde{M}$ with inclusion inducing an isomorphism $$\pi_1(C) \cong \pi_1(\tilde{M}) \cong H.$$
$C$ necessarily has nonempty boundary.  If $B \subset \partial C$ is a boundary component which is a 2-sphere, then irreducibility implies that $B$ bounds a 3-ball in $\tilde{M}$. 
That 3-ball either contains $C$ or its interior is disjoint from $C$.  The former can't happen
because that would imply the inclusion map $\pi_1(C) \to \pi_1(\tilde{M})$ is trivial.  Therefore, we can adjoin that 3-ball to $C$, removing $B$ as a boundary component and not changing 
$\pi_1(C)$. 
This process allows us to eliminate 2-spheres from  $\partial C$ and assume that 
$\partial C$ is nonempty and has infinite homology groups.  By an Euler characteristic argument, we conclude that $C$ also has infinite homology.   Then we have surjections
$H \cong \pi_1(C) \to H_1(C) \to \Z$, the required left-orderable group. \qed

A similar argument shows the following.

\begin{theorem}
Suppose $M$ is an orientable irreducible 3-manifold (possibly with boundary) such that $H_1(M)$ is infinite.  Then $\pi_1(M)$ is locally indicable.
\end{theorem}

\begin{corollary}
Knot groups are locally indicable and therefore left-orderable.
\end{corollary}

Surgery on a knot may or may not produce a 3-manifold with left-orderable fundamental group.  For example, consider the $+1$ surgery on the right-handed trefoil as indicated in Figure 1.  This means that we remove a tubular neighbourhood $N$ of the knot and attach a solid torus $S^1 \times D^2$ to the complement of $N$ in such a way that the meridian $\{*\} \times \partial D^2$ is attached to the longitudinal curve $J$ which has linking number $+1$ with the knot.  This is Dehn's original construction of the Poincar\'e homology sphere.  This manifold has fundamental group with presentation (see \cite{KL})	
$$\langle x, z | (zx)^2 = z^3 = x^5 \rangle$$

Here, $x$ and $y$ represent meridian curves indicated in the picture and $z = xy$.  This group has order 120, and cannot be left-orderable, as it clearly has torsion elements.

\begin{figure}
\centering
\includegraphics[height=60mm]{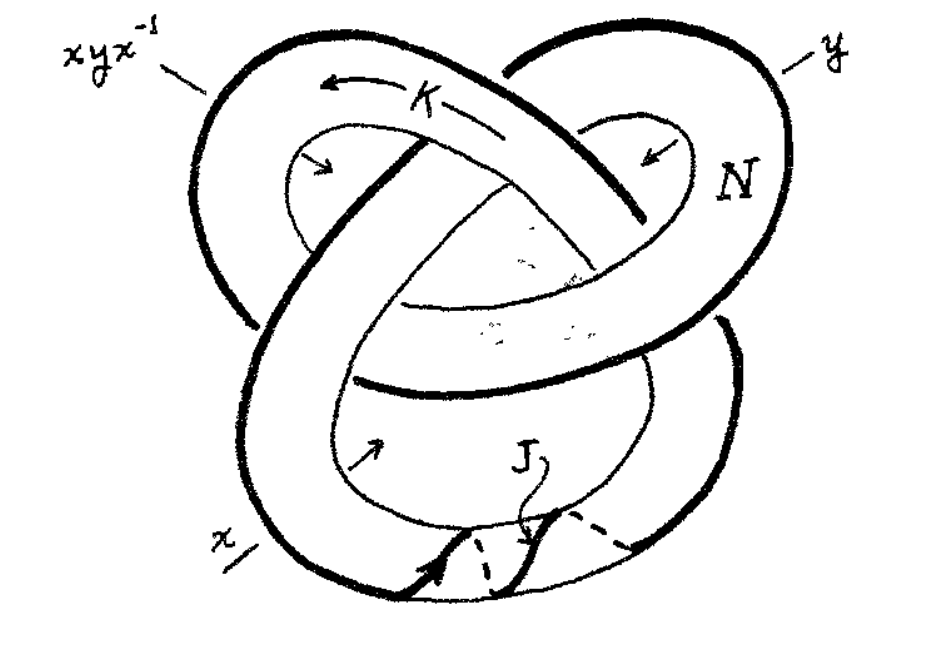}
\caption{Surgery on the trefoil}
\end{figure}

On the other hand, if we do surgery on the same knot, but along a longitudinal curve with linking number $-1$ with the knot (as $J$ in the figure, but with two more full twists at the bottom) then we get another homology sphere, with fundamental group 
$$\langle a, b | (ba)^2 = b^3 = a^7 \rangle.$$

As noted by Bergman \cite{Ber}, this group embeds in $\widetilde{SL}(2,\R)$, which we have seen is a left-orderable group .  Therefore this group is left-orderable.  Note that, since it abelianizes to the trivial group, it is not locally indicable.

\subsection{Homeomorphisms of surfaces}

Suppose $M$ is a connected triangulated surface with nonempty boundary.  Let $B$ denote the union of some or all components of the boundary, so that $B$ is nonempty.  Then define $Homeo(M,B)$ to be the group of homeomorphisms of $M$ to itself which are pointwise fixed on $B$.  The group operation is composition.  Also let $Homeo_{PL}(M,B)$ denote the subgroup of $Homeo(M,B)$ consisting of piecewise-linear homeomorphisms.

\begin{theorem}\label{homeoLO}
$Homeo_{PL}(M,B)$ is left-orderable.
\end{theorem}

\begin{figure}
\centering
{\includegraphics[height=70mm]{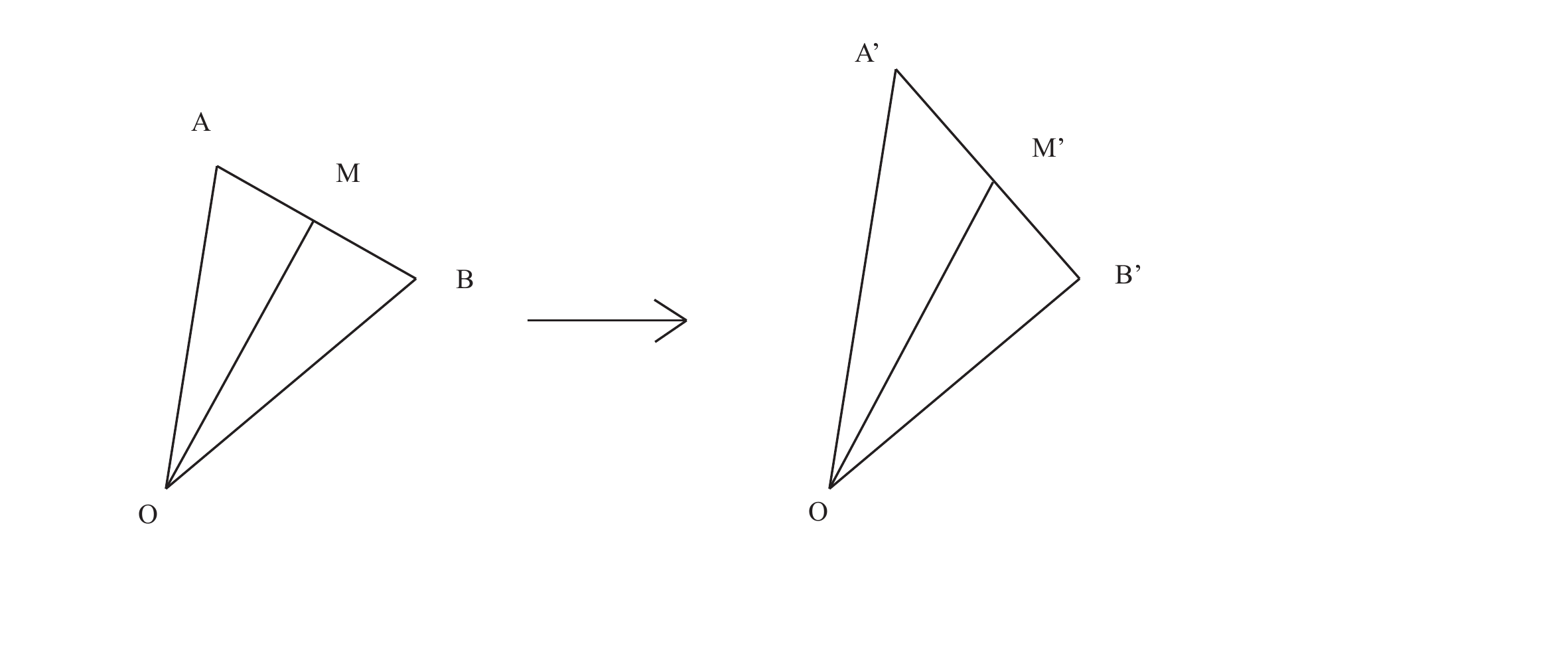}}
\caption{Showing the map of germs is nontrivial}
\end{figure}

{\bf Proof:}  We adapt an argument of Danny Calegary to this theorem for the case $Homeo_{PL}(I^2,\partial I^2)$, a result he attributes to Bert Wiest and myself.  By the Burns-Hale theorem, it suffices to consider a nontrivial finitely generated subgroup 
$H$ of $Homeo_{PL}(M,B)$ and then find a left-orderable group $L$ and nontrivial homomorphism 
$h: H \to L$.  Each of the generators of $H$ is a function that fixes some polyhedral subset of $M$ which contains $B$.  The intersection $F$ of these finitely many subsets will then be a polyhedral subset of $M$ which contains $B$; $F$ is exactly the global fixed point set of $H$.  Choose a point $O$ in the middle of an edge of a 2-simplex, on the boundary of $F$, so that small round neighbourhoods of $O$ will intersect the complement of $F$ in semidisks in some fixed triangulation of $M$.  We will consider ``germs'' of functions in $H$ at $O$ in the following sense.  Any polygonal ray $R$ emanating from $O$ is taken by the elements of $H$ to a polygonal ray $h(R)$ also originating at $O$.
If the initial segment of $R$ leaves $O$ at angle $\theta$ (measured from an edge of $F$ on which $O$ lies), $h(R)$ will be a polygonal curve whose initial segment is at angle, say, 
$\lambda_h(\theta)$.  Note that rays starting into $F$ will have their initial segments fixed.  The map $h \to (1/\pi)\lambda_h$ is a homomorphism $H \to Homeo(I, \partial I)$, a left-ordered group.

We need to check that this homomorphism is nontrivial -- that is, not every ray is mapped to a polygonal ray which starts in exactly the same direction.  If that were the case, since $O$ is on the boundary of the global fixed point set, the generators cannot preserve the length of all initial rays.  
So some generator $h$ of $H$ must send a triangle $OAB$ linearly to a triangle $OA'B'$
in which the angles $\angle AOB$ and $\angle A'OB'$ are equal and $|OA| = |OB|$, but $|OA'| \ne |OB'|$ as in Figure 2.

If $M$ is the midpoint of the side $AB$, then $M'=h(M)$ is the midpoint of $A'B'.$  We leave the reader to verify by elementary geometry that the angles $\angle A'OM'$ and $\angle M'OB'$ must be unequal.  This shows that $\lambda_h$ is not the identity mapping.  I thank Tali Pinsky for this observation. \qed 

It is also true that the group of $C^1$ homeomorphisms of the disk, fixed on the boundary, is left-orderable.  This is discussed in Calegari's blog \cite{Cal}.

\begin{proposition}
$Homeo_{PL}(I^2,\partial I^2)$  is not biorderable.
\end{proposition}

To see this, consider the PL maps $f, g: I^2 \to I^2$, where $I^2$ is regarded as the square in the $xy$-plane with $0 \le x,y\le 1$.  The map $f$ is fixed on the boundary of the square and rotates an inner square $1/4 \le x, y \le 3/4$ by 180 degrees.  The map $g$ is the identity outside the inner square, takes the point $E = (3/8, 5/8)$ to the point $E' = (5/8, 3/8)$ and extended linearly on the four triangles formed by $E$ and the sides of the inner square, taking them to the triangles formed by $E'$ and the four sides of the inner square.  One easily checks that $f^{-1}gf = g^{-1}$.  As with the Klein bottle group, the existence of a biordering on $Homeo_{PL}(I^2,\partial I^2)$ would lead to the contradiction that $g$ is positive iff $g^{-1}$ is positive. \qed

{\bf Open Question:} Is $Homeo(I^2,\partial I^2)$ left-orderable?

We note that K\'er\'ekjart\`o \cite{Ker} showed in 1920 that $Homeo(I^2,\partial I^2)$ is torsion-free.
See \cite{CK} for a discussion of this and similar results.
 
{\bf Added in proof:}  Theorem \ref{homeoLO} has been generalized in \cite{CR14} to show that
$Homeo_{PL}(M,B)$ is actually locally indicable.  There is also a version in higher dimensions.

 \bibliography{Topviewordgps}
 \bibliographystyle{plain}

\end{document}